# Time-Domain Analysis of Left-Handed Materials Based on a Dispersive Meshless Method with PML Absorbing Boundary Condition

Sheyda Shams, Farzad Mohajeri, and Masoud Movahhedi

*Abstract*—In this paper, we have proposed a dispersive formulation of scalar-based meshless method for time-domain analysis of electromagnetic wave propagation through left-handed (LH) materials. Moreover, we have incorporated Berenger's perfectly matched layer (PML) absorbing boundary condition (ABC) into the dispersive formulation to truncate open-domain structures. In general, the LH medium as a kind of dispersive media can be described by frequency-dependent constitutive parameters. The most appropriate numerical techniques for analysis of LH media are dispersive formulations of conventional numerical methods. In comparison to the conventional grid-based numerical methods, it is proved that meshless methods not only are strong tools for accurate approximation of derivatives in Maxwell's equations but also can provide more flexibility in modeling the spatial domain of problems. However, we have not seen any reports on using dispersive forms of meshless methods for simulation of wave propagation in metamaterials and applying any PML ABCs to dispersive formulation of meshless method. The proposed formulation in this paper enables us to take advantage of meshless methods in analysis of LH media. For modeling the frequency behavior of the medium in the proposed dispersive formulation, we have used auxiliary differential equation (ADE) method based on the relations between electromagnetic fields intensities and current densities. Effectiveness of the proposed formulation is verified by a numerical example; also, some basic factors which affect the accuracy and computational cost of the simulations are studied.

## I. Introduction

RESEARCH on metamaterials has been a remarkable topic for scientists in recent two decades due to their unique properties and applications. Different kinds of materials are categorized based on their constitutive parameters (permittivity and permeability) into four types: double positive (DPS), epsilon negative (EPS), mu negative (MNG), and double negative (DNG) or left-handed (LH) materials [1]. Among different kinds of materials, LH media have attracted a great attention by exhibition of extraordinary electromagnetic and optical effects. As the Poynting vector of a plane wave in LH media is antiparallel to its phase velocity vector, such medium can exhibit the unusual features of back-ward wave propagation and negative refraction.

It is very important to use proper numerical methods for accurate simulation of the electromagnetic behavior of LH media in antenna and microwave engineering. The conventional numerical methods provide precise analysis of a linear dispersive medium, with frequency-dependent constitutive parameters, only at or about the center of the desired frequency band [2]. Similarly, in [3], it is discussed that the simulation of a LH medium using the conventional finite different time-domain (FDTD) method can produce unstable results. Accordingly, the most comprehensive and reliable method for analysis of the general LH media is the dispersive formulation of conventional numerical techniques.

The dispersive formulations of mesh/grid-based numerical methods such as the dispersive FDTD methods [4] are well-developed for analysis of dispersive media including metamaterial. Beside grid-based numerical techniques, the time-domain meshless method [5] is proposed as a powerful numerical method for representing the spatial domain of nondispersive media with various boundary shapes by means of only a set of points which are called *field nodes*. Extension of meshless methods [6] for analysis of dispersive media is of a great importance due to their interesting features. Some advantages of meshless methods over FDTD methods can be summarized as [7]: in meshless methods, expansion of field variables using higher order basis functions reduces the numerical errors; these methods can provide an easier way to mesh refinements and modeling spatial domain of problems with complex geometries without mentioning connection information between nodes; meshless methods can be hybridized with other numerical methods in a simpler manner.

To benefit from the advantages of meshless methods in analysis of metamaterial, we have developed a comprehensive dispersive formulation of meshless method with an efficient absorbing boundary condition (ABC) in this paper. Prior to this study, we have seen no research on presentation of three-dimensional (3D) frequency-dependent formulation of meshless methods with any ABCs. Here, we have used Berenger's perfectly matched layer (PML) ABC [8] in our proposed dispersive meshless formulation. Development of an

ABC for truncation of the problem domain or analysis of open-region problems has a considerable importance in reducing the computational cost of the dispersive meshless method. In Berenger's PML, the problem domain is truncated using a lossy layer which is defined by electric and magnetic conductivities along each coordinate axis. Also, each component of the electric or magnetic fields is split into two orthogonal components. To satisfy the PML boundary condition, the lossy layer should be perfectly matched to the solution space [9].

The LH media can be characterized using Drude or Lorentz dispersion models. To incorporate frequency dispersion of LH media into the meshless method, we have used auxiliary differential equation (ADE) scheme. Among various techniques which are proposed for considering the frequency behavior of dispersive media in FDTD method, ADE scheme not only is an efficient technique but also is simpler than the others [10]. One of the usual approaches for implementing ADE method is based on the constitutive relation between the electric / magnetic field intensity, *i.e.,* **E** (**H**) and electric/magnetic flux density, *i.e.,* **D** (**B**). As this approach provides a second-order differential equation, discretization of the equations in time-domain requires the central-difference approximations for the first-order and second-order time-derivatives, and also a semi-implicit scheme for constant terms. Another approach is based on the relation between **E** (**H**) and electric/magnetic current density, *i.e.,* **J** (**M**), which leads to a first order differential equation. Hence, in the second approach the discretization procedure is more straightforward than the first one. To simplify the implementation of our dispersive meshless formulation with PML ABC, we have used the relation between **E** (**H**) and **J** (**M**) in this study. The basic idea for analysis of general dispersive medium using meshless methods has been presented in our previous work [11] by applying the ADE scheme based on the relation between **E** and **D**. Also, we have proposed a dispersive formulation of the vector-based meshless method in [12].

In the following sections, we will first present the details of formulation of the proposed dispersive meshless method with PML ABCs for analysis of LH media. In section III, a numerical example is employed to demonstrate the efficiency of the proposed method, and also, computational cost of the method and the effective factors on accuracy of the results are investigated. Finally, conclusions are presented in section IV.

## II. FORMULATION OF DISPERSIVE MESHLESS METHOD WITH PML ABSORBING BOUNDARY CONDITION

In this section, we present the governing equations of electromagnetic wave propagation in LH media based on the PML modification of Maxwell's equations; and meshless formulation to solve this electromagnetic problem will be obtained.

### A. PML Modification of Maxwell's Equations for Dispersive LH Media

The modified 3D time-domain Maxwell's equations based on Berenger's PML are presented in [9], [13]. First, we have rewritten these equations in frequency-domain to obtain the dispersive formulation of meshless method with PML medium as the absorbing boundary. In 3D case, splitting each component of vector fields into two subcomponents yields twelve equations as the PML modification of Maxwell's equations. The loss in the PML medium can be defined using electric and magnetic conductivities. By specifying magnetic conductivity with $\sigma_m$, the modified Faraday's law equations in frequency-domain are as the following:

$$j\omega\mu H_{yx} + \sigma_{mz} H_{yx} = -\frac{\partial}{\partial z}(E_{xy} + E_{xz}) \quad (1)$$

$$j\omega\mu H_{yz} + \sigma_{mx} H_{yz} = \frac{\partial}{\partial x}(E_{zx} + E_{zy}) \quad (2)$$

$$j\omega\mu H_{xy} + \sigma_{mz} H_{xy} = \frac{\partial}{\partial z}(E_{yx} + E_{yz}) \quad (3)$$

$$j\omega\mu H_{xz} + \sigma_{my} H_{xz} = -\frac{\partial}{\partial y}(E_{zx} + E_{zy}) \quad (4)$$

$$j\omega\mu H_{zx} + \sigma_{my} H_{zx} = \frac{\partial}{\partial y}(E_{xy} + E_{xz}) \quad (5)$$

$$j\omega\mu H_{zy} + \sigma_{mx} H_{zy} = -\frac{\partial}{\partial x}(E_{yx} + E_{yz}) \quad (6)$$

Similarly, by specifying electric conductivity with $\sigma_e$, the modified Ampere's law equations in frequency-domain are as the following:

$$j\omega\varepsilon E_{yx} + \sigma_{ez} E_{yx} = \frac{\partial}{\partial z}(H_{xy} + H_{xz}) \quad (7)$$

$$j\omega\varepsilon E_{yz} + \sigma_{ex} E_{yz} = -\frac{\partial}{\partial x}(H_{zy} + H_{zx}) \quad (8)$$

$$j\omega\varepsilon E_{xy} + \sigma_{ez} E_{xy} = -\frac{\partial}{\partial z}(H_{yx} + H_{yz}) \quad (9)$$

$$j\omega\varepsilon E_{xz} + \sigma_{ey} E_{xz} = \frac{\partial}{\partial y}(H_{zx} + H_{zy}) \quad (10)$$

$$j\omega\varepsilon E_{zx} + \sigma_{ey} E_{zx} = -\frac{\partial}{\partial y}(H_{xz} + H_{xy}) \quad (11)$$

$$j\omega\varepsilon E_{zy} + \sigma_{ex} E_{zy} = \frac{\partial}{\partial x}(H_{yx} + H_{yz}) \quad (12)$$

To characterize the frequency behavior of LH media, we have considered Drude dispersion model [ ]. In this model, the plasma frequency of the medium is defined by $\omega_p$, and $\gamma$ is the collision frequency which shows the losses of the medium. The constitutive parameters of the medium based on Dude model can be expressed as

$$\varepsilon(\omega) = \varepsilon_0 (1 - \frac{\omega_p^2}{\omega^2 - j\gamma\omega}) \quad (13)$$

$$\mu(\omega) = \mu_0 (1 - \frac{\omega_p^2}{\omega^2 - j\gamma\omega}) \quad (14)$$

Then, we substitute the permittivity and permeability of LH medium into the modified Maxwell's equations. Here, we have presented only the resulted equations for (1) and (7), as the other equations can be derived in the same manner:

$$j\omega\mu_0 H_{yx} + \sigma_{mz} H_{yx} = -\frac{\partial}{\partial z}(E_{xy} + E_{xz}) - M_{yx} \quad (15)$$

$$j\omega\varepsilon_0 E_{yx} + \sigma_{ez} E_{yx} = \frac{\partial}{\partial z}(H_{xy} + H_{xz}) - J_{yx} \quad (16)$$

where $M_{yx}$ and $J_{yx}$ are the magnetic and electric current densities which can be expressed as

$$M_{yx} = \frac{\mu_0 \omega_p^2}{(j\omega + \gamma)} H_{yx} \quad (17)$$

$$J_{yx} = \frac{\varepsilon_0 \omega_p^2}{(j\omega + \gamma)} E_{yx} \quad (18)$$

To obtain the equations in the time-domain, we calculate the inverse Fourier transform of (15)-(18). By noticing that multiplication of $j\omega$ in frequency-domain is equivalent to the time derivatives in the time-domain, the time-domain equations are obtained as:

$$\mu_0 \frac{\partial H_{yx}}{\partial t} + \sigma_{mz} H_{yx} = -\frac{\partial}{\partial z}(E_{xy} + E_{xz}) - K_{yx}(t) \quad (19)$$

$$\varepsilon_0 \frac{\partial E_{yx}}{\partial t} + \sigma_{ez} E_{yx} = \frac{\partial}{\partial z}(H_{xy} + H_{xz}) - J_{yx}(t) \quad (20)$$

$$\frac{\partial K_{yx}}{\partial t} = -\gamma K_{yx}(t) + \mu_0 \omega_p^2 H_{yx} \quad (21)$$

$$\frac{\partial J_{yx}}{\partial t} = -\gamma J_{yx}(t) + \varepsilon_0 \omega_p^2 E_{yx} \quad (22)$$

Analogously, we have developed twenty more equations according to (2)-(6) and (8)-(12).

On the other hand, by satisfying the matching condition of

$$\frac{\sigma_{e\kappa}}{\varepsilon} = \frac{\sigma_{m\kappa}}{\mu} \quad (23)$$

the characteristic wave impedance of the PML region will be the same as background medium. Therefore, the incident electromagnetic wave across the PML interface will have no reflection [8].

We will apply the radial basis function (RBF) [6] meshless method for discretizing spatial derivatives of (19)-(20). This method which is an efficient technique for solving differential form of Maxwell's equations will be discussed in the following.

*B. The RBF Meshless Method*

In meshless analysis of electromagnetic problems, we first discretize the spatial domain of the problem and its boundary using a set of electric/magnetic field nodes. The unknown electric/magnetic field ($E$ or $H$) at a point of interest can be approximated using the RBFs as follows [14]:

$$E(\mathbf{x}) = \sum_{j=1}^{N} \phi(\|\mathbf{x} - \mathbf{x}_j\|) a_j \quad (24)$$

where $\mathbf{x} = (x, y, z)$ and $\mathbf{x}_j = (x_j, y_j, z_j)$ show the location of the point of interest and the location of the node $j$, respectively. $r = \|\mathbf{x} - \mathbf{x}_j\|$ is the distance between $\mathbf{x}$ and $\mathbf{x}_j$. The radial basis function for expanding field variables is defined by $\phi$ and $a_j$ are the unknown expansion coefficients. Among different kinds of RBFs, we have applied the Gaussian function to our formulations. The Gaussian RBF can be expressed as [14]:

$$\phi_j = \phi(\|\mathbf{x} - \mathbf{x}_j\|) = e^{-\alpha r^2} \quad (25)$$

where $r = \|\mathbf{x} - \mathbf{x}_j\| = \sqrt{(x - x_j)^2 + (y - y_j)^2 + (z - z_j)^2}$ and $\alpha$ is the shape parameter of the RBF. The shape parameter plays a significant role in accuracy and stability of the meshless method's solutions. We will discuss about the impact of this parameter on the results in Section III.

The unknown expansion coefficients ($a_j$) can be obtained by supposing that (24) should be satisfied at each of $N$ nodes in the vicinity of the point of interest $\mathbf{x}$. This assumption provides a set of equations which is consisted of $N$ linear equations corresponding to all the nodes [6]. We can rewrite the compact matrix form of the obtained equations as:

$$\mathbf{A}.\mathbf{a} = \mathbf{E}_s \quad (26)$$

where the vector $\mathbf{E}_s = [E_1 \; E_2 \; ... \; E_N]^T$, $E_j$ shows the value of electromagnetic field at node $j$ ($j=1,2,...,N$), $\mathbf{a} = [a_1 \; a_2 \; ... \; a_N]^T$ contains the unknown expansion coefficients of (24), and

$$\mathbf{A} = \begin{bmatrix} \phi(\|\mathbf{x}_1 - \mathbf{x}_1\|) & \phi(\|\mathbf{x}_1 - \mathbf{x}_2\|) & ... & \phi(\|\mathbf{x}_1 - \mathbf{x}_N\|) \\ \phi(\|\mathbf{x}_2 - \mathbf{x}_1\|) & \phi(\|\mathbf{x}_2 - \mathbf{x}_2\|) & ... & \phi(\|\mathbf{x}_2 - \mathbf{x}_N\|) \\ \vdots & \vdots & & \vdots \\ \phi(\|\mathbf{x}_N - \mathbf{x}_1\|) & \phi(\|\mathbf{x}_N - \mathbf{x}_2\|) & ... & \phi(\|\mathbf{x}_N - \mathbf{x}_N\|) \end{bmatrix} \quad (27)$$

The matrix $\mathbf{A}$ is consisted of the Gaussian RBFs and is invertible. By inverting $\mathbf{A}$, (26) can be solved for obtaining the expansion coefficients $a_j$. By substitution of $a_j$ into (24), the unknown field $E(\mathbf{x})$ can be approximated as

$$E = \mathbf{B}\mathbf{A}^{-1}\mathbf{E}_s = \mathbf{\Phi}\mathbf{E}_s \quad (28)$$

where $\mathbf{B} = \mathbf{B}(\|\mathbf{x} - \mathbf{x}_j\|) = [\phi_1 \; \phi_2 \; ... \; \phi_N]$ is consisted of the radial basis functions, $\mathbf{\Phi} = [\Phi_1 \Phi_2 ... \Phi_N] = \mathbf{B}\mathbf{A}^{-1}$, and $\Phi_j = \Phi_j(\|\mathbf{x} - \mathbf{x}_j\|)$ is the shape function corresponding to each node $j$ ($j=1,2,...,N$). There is an analytical formula for calculating the partial derivatives of the shape function because the shape function is a continuous function. The first order partial derivatives of the shape function with respect to each spatial coordinate ($x$, $y$ or $z$) can be expressed as

$$\frac{\partial \mathbf{\Phi}}{\partial \kappa} = \frac{\partial \mathbf{B}(\|\mathbf{x} - \mathbf{x}_j\|)}{\partial \kappa} \mathbf{A}^{-1} \quad (29)$$

where $\kappa$ denotes $x$, $y$ or $z$. Analogously, the shape function of the magnetic field node and its derivatives can be obtained.

*C. Proposed Meshless Formulation for Analysis of LH Media*

In order to discretize 3D spatial domain of problem in meshless method, the interior medium and PML region should be represented using electric and magnetic field nodes. Due to the coupling nature between the electric and magnetic field nodes, the nodes will be arranged over the problem domain in such a way that each electric field node be surrounded by magnetic field nodes and vice versa. Moreover, all three components of the electric or magnetic field nodes will be defined at the same electric or magnetic field nodes [5].

To approximate the partial spatial derivatives of electromagnetic fields in governing equations of the problem, first, we obtain the shape functions of RBF meshless method corresponding to the field nodes through (28). Then, we calculate the spatial derivatives of electromagnetic fields using (29). On the other hand, we approximate the time derivatives in (19)-(22) (and also the obtained twenty equations corresponding to (2)-(6) and (8)-(12)) utilizing the central-difference scheme. For employing the leapfrog time-stepping scheme [9], we discretize the electromagnetic fields in time-domain by the following equations:

$$\left.\frac{\partial \mathbf{H}}{\partial t}\right|^{(q\Delta t)} \approx \frac{\mathbf{H}^{(q+\frac{1}{2})\Delta t} - \mathbf{H}^{(q-\frac{1}{2})\Delta t}}{\Delta t} \tag{30}$$

$$\left.\frac{\partial \mathbf{E}}{\partial t}\right|^{(q+\frac{1}{2})\Delta t} \approx \frac{\mathbf{E}^{(q+1)\Delta t} - \mathbf{E}^{(q)\Delta t}}{\Delta t} \tag{31}$$

Moreover, we update each component of electric and magnetic current densities in each time-step as follows:

$$\left.\frac{\partial \mathbf{M}}{\partial t}\right|^{(q+\frac{1}{2})\Delta t} \approx \frac{\mathbf{M}^{(q+1)\Delta t} - \mathbf{M}^{(q)\Delta t}}{\Delta t} \tag{32}$$

$$\left.\frac{\partial \mathbf{J}}{\partial t}\right|^{(q+1)\Delta t} \approx \frac{\mathbf{J}^{(q+\frac{3}{2})\Delta t} - \mathbf{J}^{(q+\frac{1}{2})\Delta t}}{\Delta t} \tag{33}$$

Finally, the 3D dispersive time-domain formulations of meshless method based on the Berenger's PML for analysis of LH media are as follows:

$$H_{yx,i}^{q+\frac{1}{2}} = -\frac{1}{1+\frac{\sigma_{mz}(i)\Delta t}{2\mu_0}}\left(\frac{\Delta t}{\mu_0}\right)\left(\sum_j \left(E_{xy,j}^q \partial z\Phi_j + E_{xz,j}^q \partial z\Phi_j\right)\right)$$

$$+\frac{1-\frac{\sigma_{mz}(i)\Delta t}{2\mu_0}}{1+\frac{\sigma_{mz}(i)\Delta t}{2\mu_0}} H_{yx,i}^{q-\frac{1}{2}} - \frac{1}{1+\frac{\sigma_{mz}(i)\Delta t}{2\mu_0}}\left(\frac{\Delta t}{\mu_0}\right) M_{yx,i}^q$$

$$\tag{34}$$

$$M_{yx,i}^{q+1} = \frac{1-\frac{\gamma\Delta t}{2}}{1+\frac{\gamma\Delta t}{2}} M_{yx,i}^q + \frac{\mu_0\Delta t\omega_p^2}{1+\frac{\gamma\Delta t}{2}} H_{yx,i}^{q+\frac{1}{2}} \tag{35}$$

$$H_{yz,i}^{q+\frac{1}{2}} = \frac{1}{1+\frac{\sigma_{mx}(i)\Delta t}{2\mu_0}}\left(\frac{\Delta t}{\mu_0}\right)\left(\sum_j \left(E_{zx,j}^q \partial x\Phi_j + E_{zy,j}^q \partial x\Phi_j\right)\right)$$

$$+\frac{1-\frac{\sigma_{mx}(i)\Delta t}{2\mu_0}}{1+\frac{\sigma_{mx}(i)\Delta t}{2\mu_0}} H_{yz,i}^{q-\frac{1}{2}} - \frac{1}{1+\frac{\sigma_{mx}(i)\Delta t}{2\mu_0}}\left(\frac{\Delta t}{\mu_0}\right) M_{yz,i}^q$$

$$\tag{36}$$

$$M_{yz,i}^{q+1} = \frac{1-\frac{\gamma\Delta t}{2}}{1+\frac{\gamma\Delta t}{2}} M_{yz,i}^q + \frac{\mu_0\Delta t\omega_p^2}{1+\frac{\gamma\Delta t}{2}} H_{yz,i}^{q+\frac{1}{2}} \tag{37}$$

$$H_{xy,i}^{q+\frac{1}{2}} = \frac{1}{1+\frac{\sigma_{mz}(i)\Delta t}{2\mu_0}}\left(\frac{\Delta t}{\mu_0}\right)\left(\sum_j \left(E_{yx,j}^q \partial z\Phi_j + E_{yz,j}^q \partial z\Phi_j\right)\right)$$

$$+\frac{1-\frac{\sigma_{mz}(i)\Delta t}{2\mu_0}}{1+\frac{\sigma_{mz}(i)\Delta t}{2\mu_0}} H_{xy,i}^{q-\frac{1}{2}} - \frac{1}{1+\frac{\sigma_{mz}(i)\Delta t}{2\mu_0}}\left(\frac{\Delta t}{\mu_0}\right) M_{xy,i}^q$$

$$\tag{38}$$

$$M_{xy,i}^{q+1} = \frac{1-\frac{\gamma\Delta t}{2}}{1+\frac{\gamma\Delta t}{2}} M_{xy,i}^q + \frac{\mu_0\Delta t\omega_p^2}{1+\frac{\gamma\Delta t}{2}} H_{xy,i}^{q+\frac{1}{2}} \tag{39}$$

$$H_{xz,i}^{q+\frac{1}{2}} = -\frac{1}{1+\frac{\sigma_{my}(i)\Delta t}{2\mu_0}}\left(\frac{\Delta t}{\mu_0}\right)\left(\sum_j \left(E_{zx,j}^q \partial y\Phi_j + E_{zy,j}^q \partial y\Phi_j\right)\right)$$

$$+\frac{1-\frac{\sigma_{my}(i)\Delta t}{2\mu_0}}{1+\frac{\sigma_{my}(i)\Delta t}{2\mu_0}} H_{xz,i}^{q-\frac{1}{2}} - \frac{1}{1+\frac{\sigma_{my}(i)\Delta t}{2\mu_0}}\left(\frac{\Delta t}{\mu_0}\right) M_{xz,i}^q$$

$$\tag{40}$$

$$M_{xz,i}^{q+1} = \frac{1-\frac{\gamma\Delta t}{2}}{1+\frac{\gamma\Delta t}{2}} M_{xz,i}^q + \frac{\mu_0\Delta t\omega_p^2}{1+\frac{\gamma\Delta t}{2}} H_{xz,i}^{q+\frac{1}{2}} \tag{41}$$

$$H_{zx,i}^{q+\frac{1}{2}} = \frac{1}{1+\frac{\sigma_{my}(i)\Delta t}{2\mu_0}}\left(\frac{\Delta t}{\mu_0}\right)\left(\sum_j \left(E_{xy,j}^q \partial y\Phi_j + E_{xz,j}^q \partial y\Phi_j\right)\right)$$

$$+\frac{1-\frac{\sigma_{my}(i)\Delta t}{2\mu_0}}{1+\frac{\sigma_{my}(i)\Delta t}{2\mu_0}} H_{zx,i}^{q-\frac{1}{2}} - \frac{1}{1+\frac{\sigma_{my}(i)\Delta t}{2\mu_0}}\left(\frac{\Delta t}{\mu_0}\right) M_{zx,i}^q$$

$$\tag{42}$$

$$M_{zx,i}^{q+1} = \frac{1-\frac{\gamma\Delta t}{2}}{1+\frac{\gamma\Delta t}{2}} M_{zx,i}^q + \frac{\mu_0\Delta t\omega_p^2}{1+\frac{\gamma\Delta t}{2}} H_{zx,i}^{q+\frac{1}{2}} \tag{43}$$

$$H_{zy,i}^{q+\frac{1}{2}} = -\frac{1}{1+\frac{\sigma_{mx}(i)\Delta t}{2\mu_0}}\left(\frac{\Delta t}{\mu_0}\right)\left(\sum_j \left(E_{yx,j}^q \partial x\Phi_j + E_{yz,j}^q \partial x\Phi_j\right)\right)$$

$$+\frac{1-\frac{\sigma_{mx}(i)\Delta t}{2\mu_0}}{1+\frac{\sigma_{mx}(i)\Delta t}{2\mu_0}} H_{zy,i}^{q-\frac{1}{2}} - \frac{1}{1+\frac{\sigma_{mx}(i)\Delta t}{2\mu_0}}\left(\frac{\Delta t}{\mu_0}\right) M_{zy,i}^q$$

$$\tag{44}$$

$$M_{zy,i}^{q+1} = \frac{1-\frac{\gamma\Delta t}{2}}{1+\frac{\gamma\Delta t}{2}} M_{zy,i}^q + \frac{\mu_0\Delta t\omega_p^2}{1+\frac{\gamma\Delta t}{2}} H_{zy,i}^{q+\frac{1}{2}} \tag{45}$$

Similarly, the electric fields can be updated using the following equations:

$$E_{yx,i}^{q+1} = \frac{1}{1+\frac{\sigma_{ez}(i)\Delta t}{2\varepsilon_0}}\left(\frac{\Delta t}{\varepsilon_0}\right)\left(\sum_j \left(H_{xy,j}^{q+\frac{1}{2}}\partial z\Phi_j + H_{xz,j}^{q+\frac{1}{2}}\partial z\Phi_j\right)\right)$$
$$+ \frac{1-\frac{\sigma_{ez}(i)\Delta t}{2\varepsilon_0}}{1+\frac{\sigma_{ez}(i)\Delta t}{2\varepsilon_0}} E_{yx,i}^q - \frac{1}{1+\frac{\sigma_{ez}(i)\Delta t}{2\varepsilon_0}}\left(\frac{\Delta t}{\varepsilon_0}\right) J_{yx,i}^{q+\frac{1}{2}}$$
(46)

$$J_{yx,i}^{q+\frac{3}{2}} = \frac{1-\frac{\gamma\Delta t}{2}}{1+\frac{\gamma\Delta t}{2}} J_{yx,i}^{q+\frac{1}{2}} + \frac{\mu_0 \Delta t \omega_p^2}{1+\frac{\gamma\Delta t}{2}} E_{yx,i}^{q+1}$$
(47)

$$E_{yz,i}^{q+1} = -\frac{1}{1+\frac{\sigma_{ex}(i)\Delta t}{2\varepsilon_0}}\left(\frac{\Delta t}{\varepsilon_0}\right)\left(\sum_j \left(H_{zy,j}^{q+\frac{1}{2}}\partial x\Phi_j + H_{zx,j}^{q+\frac{1}{2}}\partial x\Phi_j\right)\right)$$
$$+ \frac{1-\frac{\sigma_{ex}(i)\Delta t}{2\varepsilon_0}}{1+\frac{\sigma_{ex}(i)\Delta t}{2\varepsilon_0}} E_{yz,i}^q - \frac{1}{1+\frac{\sigma_{ex}(i)\Delta t}{2\varepsilon_0}}\left(\frac{\Delta t}{\varepsilon_0}\right) J_{yz,i}^{q+\frac{1}{2}}$$
(48)

$$J_{yz,i}^{q+\frac{3}{2}} = \frac{1-\frac{\gamma\Delta t}{2}}{1+\frac{\gamma\Delta t}{2}} J_{yz,i}^{q+\frac{1}{2}} + \frac{\mu_0 \Delta t \omega_p^2}{1+\frac{\gamma\Delta t}{2}} E_{yz,i}^{q+1}$$
(49)

$$E_{xy,i}^{q+1} = -\frac{1}{1+\frac{\sigma_{ez}(i)\Delta t}{2\varepsilon_0}}\left(\frac{\Delta t}{\varepsilon_0}\right)\left(\sum_j \left(H_{yx,j}^{q+\frac{1}{2}}\partial z\Phi_j + H_{yz,j}^{q+\frac{1}{2}}\partial z\Phi_j\right)\right)$$
$$+ \frac{1-\frac{\sigma_{ez}(i)\Delta t}{2\varepsilon_0}}{1+\frac{\sigma_{ez}(i)\Delta t}{2\varepsilon_0}} E_{xy,i}^q - \frac{1}{1+\frac{\sigma_{ez}(i)\Delta t}{2\varepsilon_0}}\left(\frac{\Delta t}{\varepsilon_0}\right) J_{xy,i}^{q+\frac{1}{2}}$$
(50)

$$J_{xy,i}^{q+\frac{3}{2}} = \frac{1-\frac{\gamma\Delta t}{2}}{1+\frac{\gamma\Delta t}{2}} J_{xy,i}^{q+\frac{1}{2}} + \frac{\mu_0 \Delta t \omega_p^2}{1+\frac{\gamma\Delta t}{2}} E_{xy,i}^{q+1}$$
(51)

$$E_{xz,i}^{q+1} = \frac{1}{1+\frac{\sigma_{ey}(i)\Delta t}{2\varepsilon_0}}\left(\frac{\Delta t}{\varepsilon_0}\right)\left(\sum_j \left(H_{zx,j}^{q+\frac{1}{2}}\partial y\Phi_j + H_{zy,j}^{q+\frac{1}{2}}\partial y\Phi_j\right)\right)$$
$$+ \frac{1-\frac{\sigma_{ey}(i)\Delta t}{2\varepsilon_0}}{1+\frac{\sigma_{ey}(i)\Delta t}{2\varepsilon_0}} E_{xz,i}^q - \frac{1}{1+\frac{\sigma_{ey}(i)\Delta t}{2\varepsilon_0}}\left(\frac{\Delta t}{\varepsilon_0}\right) J_{xz,i}^{q+\frac{1}{2}}$$
(52)

$$J_{xz,i}^{q+\frac{3}{2}} = \frac{1-\frac{\gamma\Delta t}{2}}{1+\frac{\gamma\Delta t}{2}} J_{xz,i}^{q+\frac{1}{2}} + \frac{\mu_0 \Delta t \omega_p^2}{1+\frac{\gamma\Delta t}{2}} E_{xz,i}^{q+1}$$
(53)

$$E_{zx,i}^{q+1} = -\frac{1}{1+\frac{\sigma_{ey}(i)\Delta t}{2\varepsilon_0}}\left(\frac{\Delta t}{\varepsilon_0}\right)\left(\sum_j \left(H_{xz,j}^{q+\frac{1}{2}}\partial y\Phi_j + H_{xy,j}^{q+\frac{1}{2}}\partial y\Phi_j\right)\right)$$
$$+ \frac{1-\frac{\sigma_{ey}(i)\Delta t}{2\varepsilon_0}}{1+\frac{\sigma_{ey}(i)\Delta t}{2\varepsilon_0}} E_{zx,i}^q - \frac{1}{1+\frac{\sigma_{ey}(i)\Delta t}{2\varepsilon_0}}\left(\frac{\Delta t}{\varepsilon_0}\right) J_{zx,i}^{q+\frac{1}{2}}$$
(54)

$$J_{zx,i}^{q+\frac{3}{2}} = \frac{1-\frac{\gamma\Delta t}{2}}{1+\frac{\gamma\Delta t}{2}} J_{zx,i}^{q+\frac{1}{2}} + \frac{\mu_0 \Delta t \omega_p^2}{1+\frac{\gamma\Delta t}{2}} E_{zx,i}^{q+1}$$
(55)

$$E_{zy,i}^{q+1} = \frac{1}{1+\frac{\sigma_{ex}(i)\Delta t}{2\varepsilon_0}}\left(\frac{\Delta t}{\varepsilon_0}\right)\left(\sum_j \left(H_{yx,j}^{q+\frac{1}{2}}\partial x\Phi_j + H_{yz,j}^{q+\frac{1}{2}}\partial x\Phi_j\right)\right)$$
$$+ \frac{1-\frac{\sigma_{ex}(i)\Delta t}{2\varepsilon_0}}{1+\frac{\sigma_{ex}(i)\Delta t}{2\varepsilon_0}} E_{zy,i}^q - \frac{1}{1+\frac{\sigma_{ex}(i)\Delta t}{2\varepsilon_0}}\left(\frac{\Delta t}{\varepsilon_0}\right) J_{zy,i}^{q+\frac{1}{2}}$$
(56)

$$J_{zy,i}^{q+\frac{3}{2}} = \frac{1-\frac{\gamma\Delta t}{2}}{1+\frac{\gamma\Delta t}{2}} J_{zy,i}^{q+\frac{1}{2}} + \frac{\mu_0 \Delta t \omega_p^2}{1+\frac{\gamma\Delta t}{2}} E_{zy,i}^{q+1}$$
(57)

As implementing the PML by a constant conductivity can lead to numerical reflection error at the PML surface [9], we have considered a graded profile with geometric variation for electric conductivity. The electric conductivity profile of the graded PML is considered as [9]:

$$\sigma_{e\kappa}(i) = \sigma_{\kappa,\max}\left(\frac{\rho_\kappa}{d_\kappa}\right)^n$$
(58)

where $\kappa = x, y,$ or $z$, $d_\kappa$ is the thickness of PML and $\rho_\kappa$ is the depth of node $i$ inside the PML region along each coordinate. Also, $\sigma_{\kappa,\max}$ is the maximum value of electric conductivity that can be defined using the theoretical reflection coefficient $R_{th}$ and characteristic impedance of the medium $\eta$ as

$$\sigma_{\kappa,\max} = -\frac{(n+1)\ln R_{th}}{2\eta d}$$
(59)

According to the matching condition of (23) and obtained electric conductivities in (58), corresponding magnetic conductivities will be easily found.

It is worth mentioning that (34)-(57) cover both the PML region and interior medium by choice of different values of electric and magnetic conductivities, *e.g.*, if the interior medium is free-space, it can be defined by $\sigma_{e\kappa} = \sigma_{m\kappa} = 0$.

In the general case of a 3D problem domain which is surrounded by PML media, the electric conductivities should be carefully specified at faces, edges, and corners of computational domain. Along faces the electric conductivity of PML region is not equal to interior medium's conductivity only along the axis which is normal to the interface of two media. For instance, the conditions for a PML interface which is normal to x-axis can be expressed as:

$$\sigma_{ex}^{PML} \neq \sigma_{ex}^{interior\ medium} \qquad (60)$$

$$\sigma_{ey}^{PML} = \sigma_{ey}^{interior\ medium} \qquad (61)$$

$$\sigma_{ez}^{PML} = \sigma_{ez}^{interior\ medium} \qquad (62)$$

where $\varepsilon^{PML} = \varepsilon^{interior\ medium}$, $\mu^{PML} = \mu^{interior\ medium}$ and magnetic conductivities can be obtained by matching condition of (23). On the other hand, as the faces overlap along edges and the edges overlap along corners, the electric conductivities are not equal to interior medium's conductivity along two directions at edges and along three directions at corners [9].

## III. NUMERICAL RESULTS

To verify the efficiency of our proposed formulation, we will investigate electromagnetic wave propagation through a plane metamaterial slab in this section. In subsection A, the governing equations of the problem and the simulation results are presented. The effects of nodal distributions, size of time-step, and value of the RBF's shape parameter on accuracy of the results are discussed in subsections B and C. In addition, in subsection D, the computational cost of the simulation using the proposed dispersive meshless formulation is calculated.

We have considered a 2D simulation space with dimensions $0.03m \times 0.03m$ and truncated the computational domain by the PML. A left-handed slab (with the thickness of $d = 0.01m$) is located in the middle of the domain and it is surrounded by air. A point source is located at $d/2$ above the slab and the bandwidth of excitation signal is centered at $f_0 = 30GHz$. The constitutive parameters of the slab are described using Drude dispersion model by substitution of $\omega_p = 2.666 \times 10^{11}\ rad/s$ and $\gamma = 0$ into (13) and (14). These choices for plasma and collision frequencies will provide a lossless slab with constitutive parameters of $\varepsilon_r = \mu_r = -1$ (at center frequency) and it will be matched with the free-space. To avoid excitation of the other frequency components [10], the excitation signal can be a sinusoidal wave function which has a smooth turn-on part for $m$ cycles, a constant amplitude for $n$ cycles, and again a smooth turn-off part for $m$ cycles. The period of sinusoidal function is $T_p = 1/f_0$ and it can be expressed as [4]

$$f(t) = \begin{cases} g_{on}(t)\sin(\omega_0 t) & \text{for } 0 \leq t < mT_p \\ \sin(\omega_0 t) & \text{for } mT_p \leq t < (m+n)T_p \\ g_{off}(t)\sin(\omega_0 t) & \text{for } (m+n)T_p \leq t < (m+n+m)T_p \\ 0 & \text{for } (m+n+m)T_p \leq t \end{cases} \qquad (63)$$

where switching functions $g_{on}(t)$ and $g_{off}(t)$ are the turn-on and turn-off parts of the excitation signal and given by

$$g_{on}(t) = 10x_{on}^3 - 15x_{on}^4 + 6x_{on}^5 \qquad (64)$$

$$g_{off}(t) = 1 - \left[10x_{off}^3 - 15x_{off}^4 + 6x_{off}^5\right] \qquad (65)$$

with $x_{on} = 1 - (mT_p - t)/mT_p$, $x_{off} = [t - (m+n)]/mT_p$, and by supposing $m=5$ and $n=10$.

### A. Governing Equations of the Problem

In this section, we will set the equations of a PML medium for propagation of transverse magnetic (TM) wave in a 2D problem which has no variation along $z$ axis. To discretize the spatial domain of the problem, we have used $61 \times 61$ field nodes with regular distribution in the $(x,y)$ plane. The optimum value of nodal spacing and size of time-step have been discussed in subsection B. Also, a three-layer medium is considered as the PML to truncate the computational domain and is specified by three nodes along each coordinates $x$ and $y$. Fig. 1 shows nodal distribution of field nodes over the problem domain.

For propagation of TM wave in spatial domain of the present problem, three components of electromagnetic fields are involved (i.e., $E_z$, $H_x$, and $H_y$). According to Section II, $E_z$ field should be split into two subcomponents, i.e., $E_{zx}$ and $E_{zy}$. To specify electric conductivities for inclusion of loss parameters in the PML medium; first, the maximum value of electric conductivity can be obtained by substitution of $n = 2$ and $R_{th} \sqcap 10^{-3}$ in (59). Then, the corresponding $\sigma_{ex}(i)$ and $\sigma_{ey}(i)$ can be achieved by means of (58). Finally, the dispersive meshless formulation of TM case is expressed as

$$H_{y,i}^{q+\frac{1}{2}} = \frac{1}{1+\frac{\sigma_{mx}(i)\Delta t}{2\mu_0}}\left(\frac{\Delta t}{\mu_0}\right)\left(\sum_j \left(E_{zx,j}^q \partial x\Phi_j + E_{zy,j}^q \partial x\Phi_j\right)\right)$$

$$+\frac{1-\frac{\sigma_{mx}(i)\Delta t}{2\mu_0}}{1+\frac{\sigma_{mx}(i)\Delta t}{2\mu_0}}H_{y,i}^{q-\frac{1}{2}} - \frac{1}{1+\frac{\sigma_{mx}(i)\Delta t}{2\mu_0}}\left(\frac{\Delta t}{\mu_0}\right)M_{y,i}^q$$

$$(66)$$

$$M_{y,i}^{q+1} = \frac{1-\frac{\gamma\Delta t}{2}}{1+\frac{\gamma\Delta t}{2}}K_{y,i}^q + \frac{\mu_0\Delta t\omega_p^2}{1+\frac{\gamma\Delta t}{2}}H_{y,i}^{q+\frac{1}{2}} \qquad (67)$$

$$H_{x,i}^{q+\frac{1}{2}} = \frac{1}{1+\frac{\sigma_{my}(i)\Delta t}{2\mu_0}}\left(\frac{\Delta t}{\mu_0}\right)\left(\sum_j \left(E_{zx,j}^q \partial y\Phi_j + E_{zy,j}^q \partial y\Phi_j\right)\right)$$

$$+\frac{1-\frac{\sigma_{my}(i)\Delta t}{2\mu_0}}{1+\frac{\sigma_{my}(i)\Delta t}{2\mu_0}}H_{x,i}^{q-\frac{1}{2}} - \frac{1}{1+\frac{\sigma_{my}(i)\Delta t}{2\mu_0}}\left(\frac{\Delta t}{\mu_0}\right)M_{x,i}^q$$

$$(68)$$

$$M_{x,i}^{q+1} = \frac{1-\frac{\gamma\Delta t}{2}}{1+\frac{\gamma\Delta t}{2}}M_{x,i}^q + \frac{\mu_0\Delta t\omega_p^2}{1+\frac{\gamma\Delta t}{2}}H_{x,i}^{q+\frac{1}{2}} \qquad (69)$$

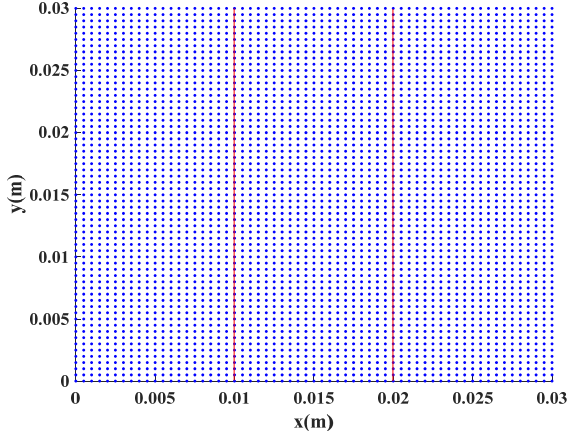

Fig. 1. Uniform nodal distribution of field nodes over the spatial domain and boundary of problem. Solid lines show the LHM slab's boundaries.

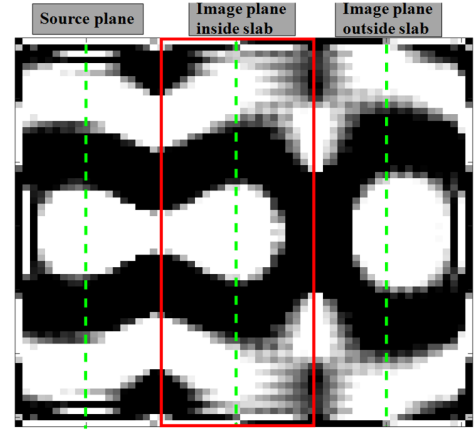

Fig. 2. The electric field intensity over meshless simulation space for a point source above a lossless LHM slab at $t = 5.7623 \times 10^{-10}$ s. The LHM slab is outlined, also source and image planes are indicated by dashed lines. In the snapshot, white (black) indicates high (low) intensity.

$$E_{zy,i}^{q+1} = -\frac{1}{1+\frac{\sigma_{ey}(i)\Delta t}{2\varepsilon_0}}\left(\frac{\Delta t}{\varepsilon_0}\right)\left(\sum_j H_{x,j}^{q+\frac{1}{2}}\partial y \Phi_j\right) + \frac{1-\frac{\sigma_{ey}(i)\Delta t}{2\varepsilon_0}}{1+\frac{\sigma_{ey}(i)\Delta t}{2\varepsilon_0}}E_{zy,i}^{q}$$

$$-\frac{1}{1+\frac{\sigma_{ey}(i)\Delta t}{2\varepsilon_0}}\left(\frac{\Delta t}{\varepsilon_0}\right)J_{zy,i}^{q+\frac{1}{2}}$$

(70)

$$J_{zy,i}^{q+\frac{3}{2}} = \frac{1-\frac{\gamma\Delta t}{2}}{1+\frac{\gamma\Delta t}{2}}J_{zy,i}^{q+\frac{1}{2}} + \frac{\mu_0\Delta t \omega_p^2}{1+\frac{\gamma\Delta t}{2}}E_{zy,i}^{q+1} \quad (71)$$

$$E_{zx,i}^{q+1} = \frac{1}{1+\frac{\sigma_{ex}(i)\Delta t}{2\varepsilon_0}}\left(\frac{\Delta t}{\varepsilon_0}\right)\left(\sum_j H_{y,j}^{q+\frac{1}{2}}\partial x \Phi_j\right) + \frac{1-\frac{\sigma_{ex}(i)\Delta t}{2\varepsilon_0}}{1+\frac{\sigma_{ex}(i)\Delta t}{2\varepsilon_0}}E_{zx,i}^{q}$$

$$-\frac{1}{1+\frac{\sigma_{ex}(i)\Delta t}{2\varepsilon_0}}\left(\frac{\Delta t}{\varepsilon_0}\right)J_{zx,i}^{q+\frac{1}{2}}$$

(72)

$$J_{zx,i}^{q+\frac{3}{2}} = \frac{1-\frac{\gamma\Delta t}{2}}{1+\frac{\gamma\Delta t}{2}}J_{zx,i}^{q+\frac{1}{2}} + \frac{\mu_0\Delta t \omega_p^2}{1+\frac{\gamma\Delta t}{2}}E_{zx,i}^{q+1} \quad (73)$$

Equations (66)-(73) would enable us to generate animated outputs for simulation of wave propagation in the structure of Fig. 1. We have analyzed the problem by means of the proposed formulation and have obtained the instantaneous distribution of $E_z$ field in Fig. 2. In our simulation, similar to the FDTD method's results for investigation of *perfect lens* effect in the LHM slab [4], two focus regions of electromagnetic wave can be found at distances $d/2$ inside the slab and $d/2$ beyond the slab. Moreover, the images do not stay stable with time and they move back and forth over time or both of them sometimes disappear. In Fig. 3, the simulated field intensity is depicted in the first and second focal planes. The capability of left-handed slab to provide near-perfect imaging of a point source can be realized by these results.

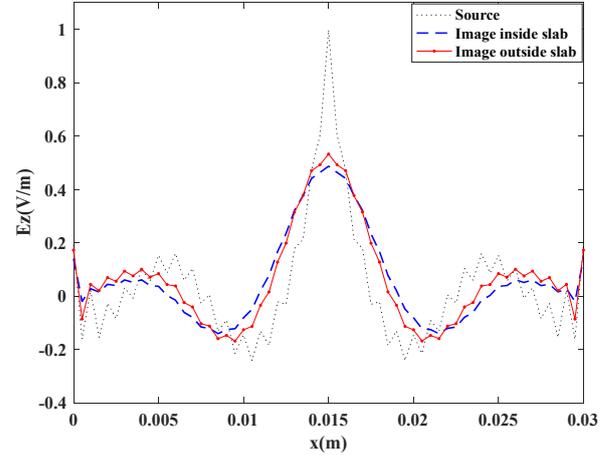

Fig. 3. Snapshots of the electric field in the source plane, the image plane inside the slab, and beyond the slab at $t = 5.7623 \times 10^{-10}$ s.

### B. The Nodal Spacing and Size of Time-Step

In the conventional RBF meshless method, the largest time-step size can be determined by the minimum of nodal spacing between any two nodes, i.e., $d_{min}$, based on the following formula [15]:

$$\Delta t \leq d_{\min}\sqrt{\varepsilon\mu} \quad (74)$$

Hence, the proper choice of nodal spacing in spatial domain discretization not only reduces numerical errors but also increases computational efficiency. In analysis of the LHM slab using the proposed dispersive meshless formulation, we examined various uniform nodal spacing to obtain the optimum nodal distance. We found the nodal spacing of $d_{\min} = \lambda_0/20$, like Fig. 1, as an optimum value which can ensure numerical accuracy and reduce the computational cost. Moreover, it provides a proper resolution for generating 2D animated outputs. On the other hand, we attained the maximum size of time-step, related to $d_{\min} = \lambda_0/20$, equal to $\Delta t = (d_{\min}/2)\sqrt{\varepsilon\mu}$.

To evaluate the accuracy of meshless simulations with different nodal distributions, we also analyzed the problem using a dispersive FDTD method and considered the results of FDTD method as the reference solution. In Fig. 4, for instance, we have illustrated that the regular nodal distribution of 61 nodes along each coordinates x and y in problem domain provides more accurate solutions than a uniform distribution of 31 nodes ($d_{min} = \lambda_0/10$). In the FDTD method, we discretized the spatial domain with a fine grid size of $\Delta x = \Delta y = \lambda_0/100$.

The $L^2$ norm of the error as a function of time is computed by

$$L^2 = \sum_{j=1}^{N} \left[ E_z(j) - E_{zr}(j) \right]^2 \tag{75}$$

In (75), $E_z(j)$ and $E_{zr}(j)$ are the normalized electric fields in $N$ points of image plane outside the slab, which are obtained by the proposed dispersive meshless formulation and the FDTD method, respectively.

*C. Impact of the Shape Parameter on RBFs*

The shape parameter, *i.e.,* $\alpha$ controls the extent of the radial scalar basis function in (25). This parameter influences the condition number of matrix **A** in (27) and also the accuracy of interpolation. Moreover, choosing smaller values for the shape parameter improves the accuracy of interpolation [16].

Due to the key role of the shape parameter in accuracy and stability of the results, we investigated the effect of choosing different values of $\alpha$ on accuracy of obtained field intensity at the image plane outside the LHM slab. We considered the parameter $\alpha_c$ related to the shape parameter $\alpha$ by the formula $\alpha_c = \alpha \cdot r^2$. For nodal distribution of Fig. 1, we found that $\alpha_c = 0.5$ is almost the smallest shape parameter which does not make matrix **A** singular. Using (75), we have computed the $L^2$ norm of error in electric field intensity for different shape parameter of RBFs. According to Fig. 5, the dispersive meshless method with $\alpha_c = 0.5$ leads to the smallest relative error among the other values of the shape parameter.

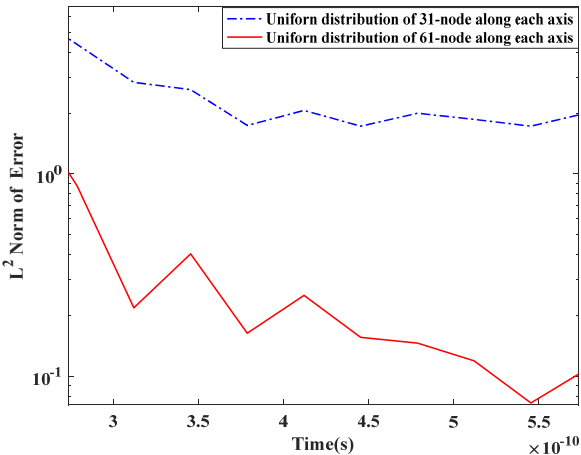

Fig. 4. The $L^2$ norm of error versus time for different nodal distributions of the proposed dispersive meshless method relative to the solution of dispersive FDTD method.

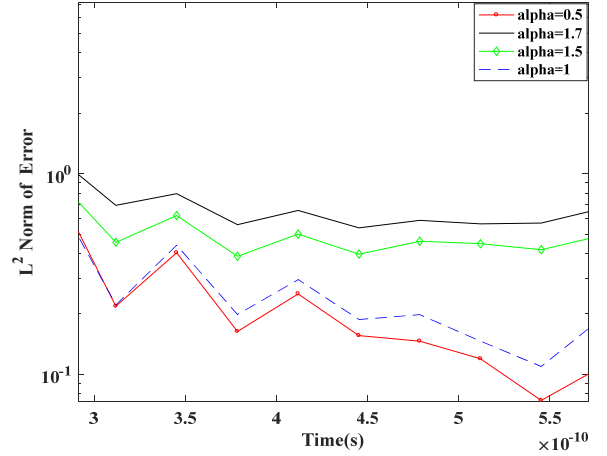

Fig. 5. The $L^2$ norm of error in electric field intensity versus time for different values of the shape parameter alpha ($\alpha_c$) in the proposed dispersive meshless formulation relative to the solution of dispersive FDTD method.

*D. Computational Cost*

To investigate the computational cost of the proposed meshless method, we considered *61×61* field nodes for the domain discretization. Then, in preprocessing stage of simulation, we obtained the shape functions of the field nodes by choosing the optimum shape parameter $\alpha_c = 0.5$. Also, the size of time-step is specified based on the proposed maximum value in subsection B. Table I lists the computational cost and memory requirement of the proposed dispersive meshless method.

TABLE I
COMPUTATIONAL COST OF THE PROPOSED METHOD

| Method | Total number of unknowns | CPU time (s) | Memory requirement (Mb) |
|---|---|---|---|
| Proposed dispersive meshless method | 3721 | 1063 | 847 |

IV. CONCLUSION

In this study, we have proposed a 3D dispersive formulation of meshless method with Berenger's PML ABC for analysis of LH materials. To simplify the incorporation of this ABC into the meshless method, the frequency behavior of LH media has been characterized by the auxiliary differential equations based on the relation between field intensities and current densities. The efficiency of the proposed formulation and also the influential parameters in accuracy and computational cost of the numerical analysis have been investigated. As a result, straightforward implementation, truncation of the problem domain by the PML, and analysis of numerical errors are some features which transform the proposed formulation into an efficient meshless technique for quantifying the characteristics of LH metamaterials and modeling electromagnetic wave interactions with them.